\documentclass{article}

\usepackage{amssymb}
\usepackage[cp1251]{inputenc}
\usepackage[english]{babel}
\usepackage{amsmath}
\usepackage{color}
\usepackage{graphics}
\usepackage{epsfig}
\usepackage[all,knot]{xy}
\xyoption{arc}

\newtheorem{proposition}{Proposition}
\newtheorem{theorem}{Theorem}

\newtheorem{corollary}{Corollary}

\def\A{{\bf A}}
\def\C{{\bf C}}
\def\F{{\bf F}}
\def\K{{\bf K}}
\def\R{{\bf R}}
\def\T{{\bf T}}

\def\cA{{\cal A}}

\def\cF{{\cal F}}
\def\cG{{\cal G}}
\def\cH{{\cal H}}
\def\cK{{\cal K}}

\def\cR{{\cal R}}
\def\cX{{\cal X}}
\def\cY{{\cal Y}}

\def\mapr#1{\smash{\mathop{\buildrel{#1}\over\longrightarrow}}}

\def\ad{{\hbox{\bf ad}}}
\def\Ad{{\hbox{\bf Ad}}}
\def\bydef{\stackrel{def}{=}}
\def\Der{{\hbox{\bf Der}\;}}
\def\Hom{{\hbox{\bf Hom}\;}}
\def\im{{\hbox{\bf Im}\;}}
\def\Int{{\hbox{\bf Int}\;}}
\def\Mor{{\hbox{\bf Mor}\;}}
\def\Obj{{\hbox{\bf Obj}\;}}
\def\Out{{\hbox{\bf Out}\;}}

\def\proof{{\bf Proof.}}
\def\qed{\hfill\vrule width2mm height2mm depth2mm}

\title{Smooth Version of Johnson's Problem Concerning Derivations of Group Algebras}
\author{A.~A.~Arutyunov and A.~S.~Mishchenko}

\begin{document}
\maketitle
\begin{abstract}
A description of the algebra of outer derivations of a group algebra of a finitely
presented discrete group is given in terms of the Cayley complex of the groupoid of
the adjoint action of the group. This task is a smooth version of Johnson's problem
concerning the derivations of a group algebra. It is shown that the algebra of outer
derivations is isomorphic to the group of the one-dimensional cohomology with compact
supports of the Cayley complex over the field of complex numbers.
\end{abstract}

\section{Introduction}

\subsection*{History of the Problem}

\subsubsection*{Inner and outer derivations}

According to numerous evidences, the derivation problem for associative algebras is
connected with B.~E.~Johnson's works concerning the cohomology of Banach algebras
(\cite{Johnson-Sinclair-1968}, \cite{Johnson-Ringrose-1969}, \cite{Johnson-1972},
\cite{Johnson-2001}). V.~Losert (\cite{Losert-2008}), who solved Johnson's problem,
formulates it in the following way as a problem concerning the derivations on group
algebras: consider the Banach algebra $\cA $and a $\cA $--bimodule
$E $. A linear mapping $$ D: \cA \mapr{}E $$ is called a derivation (or
differentiation) if, for any elements $a,b\in\cA$, the so-called Leibniz
identity (with respect to the two-sided action of the algebra $\cA$ on the
bimodule $E$) $$ D(ab)=D(a)b+aD(b), \quad a,b\in\cA. $$ (see
Definition~1.8.1 in the Dales paper (2000) \cite{Dales-2000}).

Denote the space of all derivations from the algebra $\cA$to the bimodule $E$ by $\Der(\cA,E)$. Among the derivations $\Der(\cA,E)$, we can
distinguish the so-called inner derivations $\Int(\cA,E)\subset\Der(\cA,E), $that are defined by the adjoint representations $$ \ad_{x}(a)\bydef xa-ax, \quad
x\in E, a\in\cA. $$ The derivations in the set
$\Der(\cA,E)\backslash\Int(\cA,E)$are said to be outer. It is more natural to
consider the quotient space $\Out(\cA,E)= \Der(\cA,E)/\Int(\cA,E)$as the space
of ``outer'' derivations; this space can be interpreted using the one-dimensional
Hochschild cohomology of the algebra $\cA$with coefficients in the bimodule
$E$: $$ H^{1}(\cA;E)\approx\Out(\cA,E), $$ (see the book
\cite{Pierce-1986} by R.~Pierce (1986), Definition~a, p.~248).

The derivation problem is as follows: is it true that every derivation is inner? (See
Dales(2000) \cite{Dales-2000}, (Question~5.6.B, p.~746)); i.e., is it true that $$
H^{1}(\cA;E)\approx\Out(\cA,E)=0? $$

The comparison problem for inner and outer derivations has a rich history originating
from Kaplansky's papers (1953,1958) \cite{Kaplansky-1953}, \cite{Kaplansky-1958} and
continued by Sakai (1960--1971) \cite{Sakai-1960},\cite{Sakai-1966},
\cite{Sakai-1968}, \cite{Sakai-1971}, Kadison (1966) \cite{Kadison-1966},
\cite{Kadison-Ringrose-1966}, Johnson and other authors (see, e.g.,
\cite{Johnson-Sinclair-1968}, \cite{Johnson-Ringrose-1969}, \cite{Johnson-1972}, and
\cite{Johnson-2001}).

A simpler and more natural case occurs when the bimodule $E$ is isomorphic to
the algebra $\cA$, which is certainly a bimodule over the algebra $\cA$ itself.

In the case of $E=\cA$, both the sets $\Der(\cA)=\Der(\cA,\cA)$ and
$\Int(\cA)=\Int(\cA,\cA)$ are Lie algebras with respect to the commutation
operation, and the Lie algebra $\Int(\cA)$ of inner derivations is an ideal in
the algebra $\Der(\cA)$. For this reason, it is natural to refer to the algebra
$\Out(\cA)=\Der(\cA)/\Int(\cA)$ as the algebra of outer derivations.

In this case, in Sakai's paper (1966) \cite{Sakai-1966}, it is proved that every
derivation of a $W^{*}$-algebra is inner, which is the solution of Kadison's
problem in the affirmative.

The passage to more general bimodules enabled V.~Losert to solve Johnson's problem by
finding inner derivations using elements of an algebra larger than $\cA$ rather
than elements of the original algebra $\cA.$
Namely, the derivation problem is stated as follows: are all derivations inner? This
problem was considered for the group algebras $\cA=C[G]$ of some group $G$ rather than for all algebras. More precisely, the group algebra $\overline\cA=L^{1}(G)$ of integrable functions on a locally compact group $G$ with respect to the Haar measure on $G$ together with the bimodule $E=M(G)$is considered, where $M(G)$ stands for the algebra of all bounded
measures on $G$ with the multiplication operation defined by the convolution
of measures.

In this very setting, Losert (\cite{Losert-2008}) proved that
$$
\Out(L^{1}(G),M(G))=0. $$
This result is justified by the following consideration.
For the case in which $G$ is a discrete free Abelian group with finitely many
generators, i.e., $G\approx\mathbb{Z}^{n}$, the algebra $\overline\cA=L^{1}(G)$can be identified with the Fourier algebra $A(\mathbb{T}^{n})$ of continuous functions on the $n$-dimensional torus
$\mathbb{T}^{n}$ whose Fourier coefficients form an absolutely convergent
multiple series, $\cA =A(\mathbb{T}^{n})\subset C(\mathbb{T}^{n})$(this
Fourier algebra is smaller than the algebra of all continuous functions). There are
no derivations on the algebra $A(\mathbb{T}^{n})$, since it contains
sufficiently many nonsmooth functions; however, there are no inner derivations
either, because the algebra $\overline\cA=L^{1}(G)$is commutative.

\subsection*{Choice of an appropriate class of algebras}

In the present paper, we are interested in a dense subalgebra $\cA=C[G]\subset\overline\cA$of the Banach algebra $\cA=L^{1}(G)$only
rather than in the whole algebra $\cA=L^{1}(G).$The subalgebra $\cA=C[G]$\ consists of a kind of smooth elements of the algebra $\cA=L^{1}(G).$In the choice of an appropriate class of algebras, we follow the
motivations in the paper by B.~Blackadar and J.~Cuntz (1991) \cite{Blackadar-1991} in
which the very categorial approach to the choice of the so-called smooth subalgebras
in a $C^{*}$-algebra is studied (see also the lectures of V.~Ginzburg
\cite{ginzburg-2005}, 19. Formally Smooth Algebras, p.~101)

In topology it is often important to consider an additional structure on some
topological spaces, like the smooth or piecewise linear structure. From the viewpoint
point of noncommutative geometry developed in the books of A.~Connes
\cite{Connes-1985}, \cite {Connes-1994}, an approach to the description of structures
of this kind that admits a generalization to operator algebras is to indicate some
dense *-subalgebra of the $C^*$-algebra of continuous functions. For example, a
smooth structure on a manifold $X$ can be determined by defining the subalgebra
$C^{\infty}_{0}(X)$ of the algebra $C_{0}(X)$. A piecewise-linear structure (a
triangulation) or the structure of an affine algebraic variety on $X$ can be regarded
as a choice of a suitable family of generators of the algebra $C_{0}(X)$.

When studying operator algebras, it has long been recognized that there are
circumstances in which it is natural to consider dense *-subalgebras of a given
$C^*$-algebra (in particular, in connection with cyclic cohomology or with the study
of unbounded derivations on a $C^*$-algebra.) Accepting the philosophy of
noncommutative geometry claiming that $C^*$-algebras are generalizations of
topological spaces, we can consider dense subalgebras as a tool to specify an
additional structure on the underlying space. Studies of this kind in noncommutative
differential geometry used this idea in~\cite{Connes-1985} and
\cite{Connes-Moscovici-1990}.

An important example of $C^{*}$-algebras is given by group $C^{*}$-algebras in which
subalgebras pretending to have the name of smooth subalgebras are also considered
(\cite{Jolissaint-1989}, \cite{Jolissaint-1990}, \cite{Harpe-1988}). The group
algebra $C[G]$ is an example of a dense smooth subalgebra of the group
$C^{*}$-algebra $C^{*}[G]$, and $C[G]$ is the very object of investigation in the
present paper.

\subsection*{Statement of the problem}

Consider the group algebra $\cA=C[G].$ We assume that the group $G$is a
finitely generated discrete group. Denote by $\Der(\cA)$ the set of all derivations
of the algebra~$\cA$; this set is a Lie algebra with respect to the commutator of
operators. There is a natural problem to describe all derivations of $\cA$. The inner
derivations of $\cA=C[G]$ form an ideal $\Int(\cA)\subseteq\Der(\cA)$ in the algebra
$\Der(\cA)$ of all derivations.

Johnson's problem can be generalized to the case of group algebras. Certainly,
Johnson's conjecture on the coincidence of the algebra $\Der(\cA)$ of all derivations
of the algebra $\cA=C[G]$ with the subalgebra of inner derivations fails to hold.
Therefore, Johnson's problem should be treated as the problem of evaluating the
quotient algebra of outer derivations $\Out(\A)=\Der(\cA)/\Int(\cA)$ which is
isomorphic to the one-dimensional Hochschild cohomology group of the algebra $\cA$
(with the coefficients in the bimodule $\cA$).

To every group $G$ we assign the groupoid of the adjoint action of the group~$G$,
$\cG$, and show that every derivation of the algebra $\cA=C[G]$ is uniquely defined
by an additive function on $\cG$ which satisfies some natural finiteness conditions
for the support.

For the case in which the group $G$ is finitely presented and its presentation is of
the form $G=F<X,R>$, one can transfer the presentation using the generators and
defining relations to the groupoid $\cG$: $\cG=\cF<\cX,\cR>$. This presentation
enables us to construct the Cayley complex $\cK(\cG)$ of the groupoid $\cG$ as a
two-dimensional complex whose vertices are the objects of the groupoid $\cG$, the
edges are the system of generating morphisms, and the two-dimensional cells are
formed by the system of defining relations.

Thus, the problem is to prove that the algebra of outer derivations
$\Out(\A)=\Der(\cA)/\Int(\cA)$ of the algebra $\cA$ is isomorphic to the
one-dimensional cohomology of the Cayley complex $\cK(\cG)$ of the groupoid $\cG$
with finite supports: $$ \Out(C[G])\approx H^{1}_{f}(\cK(\cG); \R) $$
(Corollary~\ref{c1} to Theorem~\ref{t4}).

\subsection*{Plan of the paper}
In Sec.~2 we consider the necessary properties of the groupoid $\cG$ of the adjoint
action of the group $G$. In Sec.~3 we establish the key theorem on the description of
derivations using characters on the groupoid $\cG$. In Sec.~4, a generalization of
the Cayley complex to the case of a presentation of the groupoid $\cG$ is constructed
and the theorem on an isomorphism between the algebra of outer derivations of the
group algebra $C[G]$ and the one-dimensional cohomology of the Cayley complex of the
groupoid $\cG$ is proved. Finally, in Sec.~5, the simplest examples of groups are
considered for which the Cayley complexes of the corresponding groupoids can be
described.

\section{Groupoid of the adjoint action of a group}
\subsection{Linear operators on the group algebra}

Consider the group algebra $\cA=C[G].$ We assume that $G$is a finitely
presented discrete group.

An arbitrary element $u\in\cA$ is a finite linear combination $u=\sum\limits_{g\in
G}\lambda^{g}\cdot g. $ Consider an arbitrary linear operator on the group algebra
$\cA=C[G]$,
$$
X:\cA\mapr{}\cA.
$$
Since the algebra $\cA$ is a vector space in which the set of all elements of the
group $G$ is a basis, it follows that the linear operator $X$ has the following
matrix form: \begin{equation}\label{1} X(u)=\sum\limits_{h\in G}\left(\sum\limits_{g\in
G}x_{g}^{h}\lambda^{g}\right)\cdot h. \end{equation} The matrix entries $x^{h}_{g}$are
given by the equations \begin{equation}\label{2} X(g)=\sum\limits_{h}x^{h}_{g}\cdot h\in\cA. \end{equation} Since
the sum in equation (\ref{2}) must be finite, this means that the matrix $X=\|x^{h}_{g}\|_{g,h\in G}$must satisfy the natural finiteness condition:
\begin{itemize}
\item[(\F1)] For every subscript $g\in G$, the set of the superscripts $h\in G$ for
which $x^{h}_{g}$ is nonzero is finite.
\end{itemize}
In particular, it follows from condition (\F1) that the outer sum in the matrix
representation (\ref{1}) is also finite.

Certainly, the converse assertion also holds: if a matrix $X=\|x^{g}_{h}\|_{g,h\in
G}$ satisfies condition (\F1), then it well defines a linear operator $X:\cA\mapr{}\cA$by formula (\ref{1}). All this justifies that both the operator
$X$and its matrix $X=\|x^{h}_{g}\|_{g,h\in G}$are denoted by the same
symbol $X.$
Consider now a so-called differentiation (derivation) in the algebra $\cA$, i.e., an
operator $X$ for which the following condition holds:
\begin{itemize}
\item[(\F2)]$ X(u\cdot v)=X(u)\cdot v+u\cdot X(v), u,v\in \cA. $
\end{itemize}

The set of all derivations of the algebra $\cA$ is denoted by $\Der(\cA)$ and forms a
Lie algebra with respect to the commutator of operators.

There is a natural problem to describe all derivations of the algebra $\cA$. To this
end, it is necessary to satisfy two conditions, (\F1) and (\F2). It is more or less
simple to verify each of the conditions separately. The simultaneous validity of
these conditions is one of the tasks of this paper.

There is a class of the so-called inner derivations, i.e., operators of the form
$$
X=\ad(u), \quad X(v)=\ad(u)(v)=[u,v]=u\cdot v-v\cdot u, \quad  u,v\in \cA.
$$
All inner derivations satisfy automatically both the conditions (\F1) and (\F2). We
denote their set by $\Int(\cA)$; this is a Lie subalgebra of the Lie algebra
$\Der(\cA)$, $$ \Int(\cA)\subseteq\Der(\cA). $$

\begin{proposition}
The subalgebra $\Int(\cA)\subseteq\Der(\cA)$is an ideal.
\end{proposition}

Indeed, we are to verify the validity of the condition $$
[\Int(\cA),\Der(\cA)]\subset \Int(\cA). $$ If $\ad(u)\in\Int(\cA)$,
$X\in\Der(\cA)$, then the commutator $[\ad(u),X]$ is evaluated by the formula $$
\begin{array}{l}
[\ad(u),X](v)=\ad(u)(X(v))-X(\ad(u)(v))=[u,X(v)]-X([u,v])=
\\=[u,X(v)]-[X(u),v]-[u,X(v)]=-\ad(X(u))(v),
\end{array}
$$ i.e., $[\ad(u),X]\in\Int(\cA)$.

Thus, the quotient space $\Out(\cA)=\Der(\cA)/\Int(\cA)$ is a Lie algebra, which is
called the algebra of outer derivations.

\subsection{Definition of the groupoid $\cG$ of the adjoint action of the group $G$}

Denote by $\cG$ the groupoid associated with the adjoint action of the group $G$ (or
the groupoid of adjoint action, see, for example, Ershov (2012) \cite{Ershov-2012},
p.~18, Example~j). The groupoid $\cG$ consists of the objects $\Obj(\cG)=G$ and the
morphisms $$ \Mor(a,b)=\{g\in G: ga=bg \hbox{ or } b=\Ad(g)(a)\}, \quad a,b\in
\Obj(\cG). $$ It is convenient to denote the elements of the set of all morphisms
$\Mor(\cG)=\coprod\limits_{a,b\in \Obj(\cG)}\Mor(a,b)$ in the form of columns $$
\xi=\left(\frac{a\mapr{}b}{g}\right)\in \Mor(a,b), \quad b=gag^{-1}=\Ad(g)(a). $$
The composition $*$ of two morphisms is defined by the formula $$
\begin{array}{l}
\left(\frac{a\mapr{}c}{g_{2}g_{1}}\right)= \left(\frac{b\mapr{}c}{g_{2}}\right)*
\left(\frac{a\mapr{}b}{g_{1}}\right),\\\\
b=\Ad(g_{1})(a), \\\\
c=\Ad(g_{2})(b)=\Ad(g_{2})(\Ad(g_{1})(a))=\Ad(g_{2}\Ad(g_{1})(a))
\end{array}
$$ which corresponds to the diagram $$ \xymatrix{
& \Ad(g_{1})(a)\ar@{=}[d]& \Ad(g_{2}g_{1})(a)\ar@{=}[d]\\
a\ar[r]^{g_{1}}\ar@/_15pt/[rr]_{g_{2}g_{1}}& b\ar[r]^{g_{2}}& c } $$

Note that the groupoid $\cG$ is decomposed into the disjoint union of its
subgroupoids $\cG_{\langle g\rangle}$ that are indexed by the conjugacy classes
$\langle g\rangle$ of the group $G$:
$$
\cG=\coprod\limits_{\langle g\rangle\in\langle G\rangle}\cG_{\langle g\rangle},
$$
where $\langle G\rangle$ stands for the set of conjugacy classes of the group~$G$.
The subgroupoid $\cG_{\langle g\rangle}$ consists of the objects $\Obj(\cG_{\langle
g\rangle})=\langle g\rangle$ and the morphisms $\Mor(\cG_{\langle
g\rangle})=\coprod\limits_{a,b\in\langle g\rangle}\Mor(a,b)$.

\subsection{Linear operators as functions on the groupoid $\cG$}

A linear operator $X:\cA\mapr{}\cA$ is described by the matrix
$X=\|x^{h}_{g}\|_{g,h\in G}$ satisfying condition (\F1). The same matrix $X$ defines
a function on the groupoid $\cG$:
$$
T^{X}:\Mor(\cG)\mapr{}R,
$$
associated with~$X$, which is defined by the following formula: if $\xi$ is a
morphism, $$ \xi=\left(\frac{a\mapr{}b}{g}\right)\in\Mor(\cG), $$ then we set
$$ T^{X}(\xi)=T^{X}\left(\frac{a\mapr{}b}{g}\right)=x^{ga=bg}_{g}. $$
Condition~(\F1) imposed on the coefficients of the matrix $X$can be
reformulated in terms of the function $T$on the morphisms $\Mor(\cG)$ of the
groupoid $\cG$ as follows:
\begin{itemize}
\item[(\T1)] for every element $g\in G$, the set of morphisms of the form $$
\xi=\left(\frac{a\mapr{}b}{g}\right) $$ for which $T^{X}(\xi)\neq 0,$is
finite.
\end{itemize}
The set of all morphisms $\Mor(\cG)$can be represented in the form of the
disjoint union of the sets $$ \Mor(\cG)=\coprod\limits_{g\in G}\cH_{g}, $$ where
$$ \cH_{g}=\left\{\xi=\left(\frac{a\mapr{}b}{g}\right):a\in G,b=gag^{-1}\in
G\right\}. $$ Then the condition (\T1) imposed on the function $T$can
equivalently be reformulated as follows:
\begin{proposition}
A function $$ T:\Mor(\cG)\mapr{}\C, $$ is defined by a linear operator $$
X:\cA\mapr{}\cA, \quad T=T^{X}, $$ if and only if, for any element $g\in
G,$\ the restriction ${\left(T\right)_{|}}_{\cH_{g}}:\cH_{g}\mapr{}\C$is a
finitely supported function (on the set $\cH_{g}$).
\end{proposition}
 We say that the functions $T:\Mor(\cG)\mapr{}\C$of this kind satisfying
the condition of finite support on every subset of the form $\cH_{g},$$g\in G,$are locally finitely supported functions on the groupoid $\cG$.
Denote the set of locally finitely supported functions on the groupoid $\cG$ by
$C_{f}(\cG)$. All this means that the correspondence $T$ assigning to every operator
$X\in\Hom(\cA,\cA)$ the function $T^{X}$ is an isomorphism between the spaces
$\Hom(\cA,\cA)$ and $C_{f}(\cG)$, i.e., the following assertion holds.
\begin{theorem}
The homomorphism
$$
T:\Hom(\cA,\cA)\mapr{}C_{f}(\cG)
$$
is an isomorphism.
\end{theorem}

\section{Derivations of the group algebra as characters on the groupoid}
The algebra of derivations $\Der(\cA)$ treated as linear operators is a subspace of
$\Hom(\cA,\cA)$. Thus, the correspondence $T$ takes the algebra of derivations
$\Der(\cA)$ to some subspace $\T_{f}(\cG)\subset C_{f}(\cG)$.

Consider two morphisms $\xi=\left(\frac{a\mapr{}b}{g_{1}}\right)$ and
$\eta=\left(\frac{b\mapr{}c}{g_{2}}\right)$, which thus admit the composition
$$
\eta*\xi=\left(\frac{a\mapr{}c}{g_{2}g_{1}}\right).
$$
\begin{theorem}
An operator $X:\cA\mapr{}\cA$ is a differentiation (i.e., a derivation) if and only
if the function $T^{X}$ (on the groupoid $\cG$) associated with the operator $X$
satisfies the additivity condition
\begin{itemize}
\item[(\T2)] $$T^{X}(\eta*\xi)=T^{X}(\eta)+T^{X}(\xi)$$
\end{itemize}
for every pair of morphisms $\xi$ and $\eta$ admitting the composition $\eta*\xi$.
\end{theorem}
\proof \ Let the matrix of the operator $X$have the form $X=\|x^{h}_{g}\|_{g,h\in G},$and thus the function $T^{X}$takes the
following value on the element $\xi$:
\begin{equation}\label{e3}
T^{X}(\xi)=T^{X}\left(\frac{a\mapr{}b}{g}\right)=x^{ga=bg}_{g}.
\end{equation}

Consider two morphisms $\xi=\left(\frac{a\mapr{}b}{g_{1}}\right),$$\eta=\left(\frac{b\mapr{}c}{g_{2}}\right)$ admitting the composition $\eta*\xi= \left(\frac{a\mapr{}c}{g_{2}g_{1}}\right).$ Then $$
T^{X}(\eta*\xi)=x_{g_{2}g_{1}}^{g_{2}g_{1}a=cg_{2}g_{1}}=x_{g_{2}g_{1}}^{h}, $$
$$ T^{X}(\xi)=x_{g_{1}}^{g_{1}a=bg_{1}}=x_{g_{1}}^{g_{2}^{-1}h}, $$ $$
T^{X}(\eta)=x_{g_{2}}^{g_{2}b=cg_{2}}=x_{g_{2}}^{hg_{1}^{-1}}. $$ On the other
hand, $$ X(g_{2}g_{1})=X(g_{2})g_{1}+g_{2}X(g_{1}). $$ In other words, $$
\begin{array}{l}
X(g_{2}g_{1})= \sum\limits_{h\in G} x^{h}_{g_{2}g_{1}}\cdot h= \sum\limits_{h\in
G}x^{h}_{g_{2}}\cdot h\cdot g_{1}+ g_{2}\cdot\sum\limits_{h\in G}x^{h}_{g_{1}}\cdot
h\\= \sum\limits_{h\in G}x^{hg_{1}^{-1}}_{g_{2}}\cdot h+ \sum\limits_{h\in
G}x^{g_{2}^{-1}h}_{g_{1}}\cdot h.
\end{array}
$$ Thus, $$ x^{h}_{g_{2}g_{1}}= x^{hg_{1}^{-1}}_{g_{2}}+x^{g_{2}^{-1}h}_{g_{1}}
$$ Finally, $$ T^{X}(\eta*\xi)=T^{X}(\eta)+T^{X}(\xi).\qed $$

Every function $T:\Mor(\cG)\mapr{}R$on the groupoid $\cG$satisfying
the additivity condition (T2) is called a {\it character}. Denote the set of all
characters on the groupoid $\cG$by $\T(\cG).$Denote the space of all
locally finitely supported characters of the groupoid $\cG$by
$\T_{f}(\cG)\subset\T(\cG).$
Thus, the correspondence $T$ defines a mapping from the algebra of derivations
$\Der(\cA)$ to the space $\T_{f}(\cG)$ of locally finitely supported characters on
the groupoid $\cG$:
\begin{theorem}\label{t3}
The mapping $$ T:\Der(\cA)\mapr{}\T_{f}(\cG), $$ is an isomorphism.
\end{theorem}

\section{Cayley complex of a groupoid}

Here we intend to apply the so-called geometric methods of combinatorial group theory
to study the problem of describing the derivations of the group algebra of a finitely
presentable discrete group. Following, for example, the book of R.~Lyndon and
P.~Schupp (1980, \cite {lyndon-schupp-1980}), one can assign to every discrete
finitely presentable group the so-called Cayley graph and its two-dimensional
generalization, the Cayley complex, which consists of the elements of the group as
vertices, of the system of generators as edges, and of the system of defining
relations as two-dimensional cells. The topological properties of the Cayley complex
are responsible for certain algebraic properties of the group $G$itself.

The geometric construction of the Cayley complex for a finitely presentable group
$G $can be generalized to the case of groupoids; in particular, to the case
of the groupoid $\cG $of the adjoint action of the group $G $. Since
the derivations of the group algebra $\Der(C[G])$can be described as
characters on the groupoid $\cG,$it follows that the topological properties
of the Cayley complex $\cK(\cG) $of the groupoid $\cG $enable us to
describe some properties of derivations.

\subsection{Presentation of a finitely presentable group}
Consider a finitely presentable group $G$, $$ G=F<X,R>, $$ where $X=\{x_{1},x_{2},\dots,x_{n}\}$is a finite set of generators and $R=\{r_{1},r_{2},\dots,r_{m}\}$is a finite set of defining relations.

By analogy with a free group (see, e.g., Kargapolov, Merzlyakov,
\cite{kargapolov-merzliakov-1982}, pp.~122--124), an arbitrary element $g\in
G$can be represented as a word $s\in S(Y),$$g=g(s)\in G,$composed
of letters of the alphabet $Y=X\sqcup X^{-1},$where $$
X^{-1}=\{x_{1}^{-1},x_{2}^{-1},\dots,x_{n}^{-1}\}, $$ i.e., $$
s=y_{1}y_{2}y_{3}\dots y_{l}, \quad y_{j}\in Y. $$

The words $s\in S(Y)$representing the same element $g=g(s)\in G$are
obtained from one another by successive operations of reduction of words and
operations inverse to reduction. Every reduction operation is as follows. Let a word
$s\in S(Y)$be represented as a concatenation of three subwords $s=s_{1}\rho s_{2},$where the middle word $\rho$is equal to one of the
following words: $$
\begin{array}{ll}
\rho=\sigma\sigma^{-1}, & \sigma\in S(Y),\\
\rho= r, & r\in R\sqcup R^{-1}.
\end{array}
$$ Here, if the word $\sigma$is of the form
$\sigma=\{y_{1}y_{2}y_{3}\dots y_{k}\},$then, by definition,
$\sigma^{-1}=\{y_{k}^{-1}\dots y_{3}^{-1}y_{2}^{-1}y_{1}^{-1}\}.$In this case,
the new word $s'=s_{1}s_{2}$is, by definition, the result of reduction of the
word $s.$After finitely many reductions, the word becomes irreducible. The
inverse operation $s'=s_{1}s_{2}\Rightarrow s_{1}\rho s_{2}$is the operation
of insertion. Two words $s_{1}$and $s_{2}$are said to be equivalent
if there is a finite sequence of operations of reduction and insertion taking one
word, $s_{1},$to another, $s_{2}. $
It is certainly necessary to prove that diverse sequences of contractions lead to
equivalent results in the form of an irreducible word. This problem, the so-called
word problem, is not always decidable \cite{lyndon-schupp-1980}.

\subsection{Presentation of the groupoid of the adjoint action of a group}

The groupoid $\cG,$whose set of objects is $\Obj(\cG)=G$and the set of
morphisms $\Mor(\cG)$consists of $$
\begin{array}{l}
\Mor(\cG)=\coprod\limits_{a,b\in\Obj(\cG)}\Mor(a,b), \\\\
\Mor(a,b)=\{g\in G:b=a^{g}=gag^{-1}\}.
\end{array}
$$ is decomposed into a disjoint sum of subgroupoids generated by the conjugacy
classes. More precisely, denote by $\langle g\rangle$the conjugacy class of an
element $g\in G$, $$ \langle g\rangle=\{g^{h}: h\in G\}. $$ Denote the set
of conjugacy classes by $\langle G\rangle$, $\langle G\rangle=\{\langle
g\rangle:g\in G\}$. The group $G$is decomposed into the disjoint union of
the conjugacy classes $$ G=\coprod\limits_{g\in G}\langle g\rangle=
\coprod\limits_{\langle g\rangle\in \langle G\rangle}\langle g\rangle. $$

Correspondingly, the groupoid $\cG$can also be decomposed into a disjoint
union of subgroupoids $$ \cG=\coprod\limits_{\langle g\rangle\in\langle
G\rangle}\cG_{\langle g\rangle} $$ that are defined by their objects and morphisms
as follows: $$
\begin{array}{l}
\Obj(\cG_{\langle g\rangle})=\langle g\rangle,\\\\
\Mor(\cG_{\langle g\rangle})=\coprod\limits_{a,b\in\langle g\rangle}\Mor(a,b).
\end{array}
$$

The finite set of generators $X=\{x_{1},x_{2},\dots,x_{n}\}$and the finite
set of defining relations $R=\{r_{1},r_{2},\dots,r_{m}\}$are transferred to
the generators and relations of the groupoid $\cG$, which we denote by
$\cX$and $\cR$. Thus, the set of morphisms $\Mor(\cG)$can be
denoted by $\cF<\cX,\cR>,$$$ \Mor(\cG)=\cF<\cX,\cR>. $$ Let us define
$\cX$as the set of all morphisms of the form $$
\cX=\left\{\xi=\left(\frac{a\mapr{}b}{x}\right): x\in X, a\in\Obj(\cG),
b=a^{x}\right\}. $$ Let $\cY=\cX\sqcup\cX^{-1};$consider $\cY$as an
alphabet, $$ \cY=\left\{\xi=\left(\frac{a\mapr{}b}{y}\right): y\in Y=X\sqcup
X^{-1}, a\in\Obj(\cG), b=a^{y}\right\}. $$ The set $S(\cY)$is the set of
all admissible words $s$in the alphabet $\cY,$i.e., words formed by
the letters of the alphabet $\cY,$$s=\xi_{1}\xi_{2}\xi_{3}\cdots\xi_{l}
$such that $$ \xi_{i}=\left(\frac{a_{i}\mapr{}a_{i+1}}{y_{i}}\right),
\quad\xi_{i}\in\cY,\quad 1\leq i\leq l. $$ Every admissible word $s\in
S(\cY)$defines a morphism $\xi(s)\in\Mor(\cG)$by the formula $$
\xi(s)=\xi_{1}*\xi_{2}*\xi_{3}*\cdots*\xi_{l}. $$ This representation of the
morphism $\xi$in the form of an admissible word $s$is not unique, and
enables one to make a reduction of the word $s$by the following rule. Define
first the system of relations $\cR$generated by the set $R$of defining
relations for the group $G.$Every relation $r_{i}\in R$is written out
in the form of a word $$ r_{i}=y_{i1}y_{i2}y_{i3}\cdots y_{il_{i}},\quad y_{ij}\in
Y. $$ The relations $r_{i}$generate the system of admissible words
$\rho_{i,a}, a\in\Obj(\cG),$of the form  $$
\begin{array}{l}
\rho_{i,a}= \left(\frac{a_{1}\mapr{}a_{2}}{y_{i1}}\right)
\left(\frac{a_{2}\mapr{}a_{3}}{y_{i2}}\right)
\left(\frac{a_{3}\mapr{}a_{4}}{y_{i3}}\right)\cdots
\left(\frac{a_{l_{i}}\mapr{}a_{1}}{y_{il_{i}}}\right), \\\\
a=a_{1}, \quad a_{j+1}=a_{j}^{y_{ij}}, \quad 1\leq j\leq l_{i}, \quad
a_{l_{i}+1}=a_{1},
\end{array}
$$ which serve as the defining relations of the groupoid $\cG.$Denote the
set of all admissible words of the form $\rho_{i,a}$by $\cR,$$$
\cR=\{\rho_{i,a}: 1\leq i\leq l_{i},\quad a\in\Obj(\cG)\}, $$ $\cR\subset S(\cY)$.
Thus, the operation of reduction is carried out for an admissible word $s$as
follows. Let an admissible word $s$be representable in the form of the
concatenation of three words $$ s=s_{1}\eta s_{2}, $$ were the middle word is
equal to one of the following words: $$
\begin{array}{ll}
\eta= \sigma\sigma^{-1}, & \sigma\in S(\cY) \\
\eta= \rho, & \rho\in \cR\sqcup\cR^{-1}\subset S(\cY).
\end{array}
$$ In this case, the result of reduction is the word $s'=s_{1}s_{2},$\
which is certainly admissible. The inverse operation $s'=s_{1}s_{2}\Rightarrow
s=s_{1}\eta s_{2}$is called the operation of admissible insertion.

Thus, two admissible words $s$and $s'$define the same morphism, i.e.,
$$ \xi(s)=\xi(s')\in\Mor(\cG), $$ if and only if the words are equivalent, $s\sim s'$, i.e., when there is a finite sequence of operations of two types:

1) the operation of reduction,

2) the operation of admissible insertion.

\subsection{Construction~of~the~Cayley~complex~of~the~groupoid~$\cG$}

\subsubsection{Cayley complex of a group $G$}
Before constructing the Cayley complex of the groupoid $\cG$by analogy with
the Cayley complex of the group $G$itself, recall the construction of the
Cayley complex of $G$ from its presentation in the form of finitely many
generators $X$and finitely many defining relations $R$, $\cF(X)/R$. We
follow the book by Lyndon and Schupp (1980, \cite{lyndon-schupp-1980}, p.~174,
Chap.~3, \S~4, Cayley complexes). The group $G$is treated there as the
groupoid of the action of $G$on itself with the help of right multiplication:
the action of an element $g\in G $on $G $is given by the rule
 $$ G\times G\mapr{}G,
\forall g\in G, h\mapsto hg,\quad h\in G. $$ Thus, we obtain the groupoid of this
action, say, $r\cG$. The Cayley complex is constructed in the book by Lyndon
and Schupp from this very groupoid $r\cG$rather than from the group $G$. The
objects of this category are the elements $h\in G$of $G$themselves,
$\Obj(r\cG)=G$, and the morphisms $\Mor(r\cG)$are right shifts on the
group $G$, $h\mapsto hg$. This means that the set $\Mor(a,b)$consists of precisely one element $g\in G$, namely, $g=a^{-1}b.$It can
readily be seen that the category $r\cG$thus constructed is a groupoid. Thus,
the Cayley complex of the group $G$defined in the book by Lyndon and Schupp
(1980, \cite{lyndon-schupp-1980}, p.~174) is in fact constructed from the groupoid
$r\cG,$and the construction by itself can be generalized to arbitrary
groupoids associated with an action of the group $G$.

By the definition in the book by Lyndon and Schupp (1980, \cite{lyndon-schupp-1980},
p.~174), for the groupoid $r\cG$of the right action of the group $G$, the
Cayley complex $\cK(r\cG)$consists of vertices, edges, and two-dimensional
cells. The set of vertices $\cK_{0}(r\cG)$is the set of all objects of the
groupoid $r\cG,$$\cK_{0}(r\cG)=\Obj(r\cG)\approx G.$
The set of edges of the groupoid $r\cG$, $\cK_{1}(r\cG)$, is formed by the morphisms
of the form $\xi=\left(\frac{a\mapr{}ag}{g}\right),$ $g\in X\sqcup X^{-1},$ i.e.,
$\xi\in\cX\sqcup\cX^{-1}=\cY.$ The edges $\xi=\left(\frac{a\mapr{}ag}{g}\right)$ and
$\eta=\left(\frac{ag\mapr{}a}{g^{-1}}\right)$ are assumed to be the same edge with
opposite orientation. Thus, the edges $\xi\in\cK_{1}(r\cG)$ are defined by the set
$\cX$ of generators of the groupoid $r\cG$.

The set of two-dimensional cells, $\cK_{2}(r\cG)$, is defined using sequences of
morphisms defined by words $\rho\in\cR\sqcup\cR^{-1}\subset
S(\cX\sqcup\cX^{-1})=S(\cY)$. The two-dimensional cells are the planar orientable
polygons $\sigma(\rho)$defined by the words $\rho\in\cR\sqcup\cR^{-1}$that determine the boundaries of the polygons $\sigma(\rho)$as closed cycles
formed by the edges of the word $\rho.$The cells $\sigma(\rho)$and
$\sigma(\rho^{-1})$are assumed to be equal and have opposite orientation. The
two-dimensional cells $\sigma(\rho)$are pasted to the 1-skeleton of the
complex $\cK(\cG_{[c]}$by the natural identification of the edges of the
boundary of the cell $\sigma(\rho)$with the corresponding edges of the complex
$\cK(\cG_{[c]}),$preserving the orientation.

\subsubsection{Cayley complex of the groupoid $\cG$of adjoint action}
The only difference between the groupoid $\cG$ and the groupoid $r\cG$ is that the
former is defined by another action of the group $G$, namely, the adjoint action:
$\Ad_{g}(a)=dad^{-1}$, $g,h\in G$. Therefore, the Cayley complex of the groupoid
$\cG$ is constructed by analogy with the Cayley complex of the group $G$.

Namely, the vertices, i.e., the zero-dimensional cells $\cK_{0}(\cG)$of the
complex $\cK(\cG),$are the objects, $a\in\Obj(\cG)\approx G.$
The one-dimensional edges, i.e., the oriented cells of dimension 1, $\cK_{1}(\cG)$,
joining vertices $a$and $b$, are the morphisms $\xi\in\Mor(a,b)$of the form  $$ \xi=\left(\frac{a\mapr{}b}{y}\right), \quad y\in Y=X\sqcup X^{-1},
\quad a\in\Obj(\cG), \quad b=yay^{-1}\in\Obj(\cG). $$ The set of edges described
above is denoted by $\cX$; let $\cY=\cX\sqcup\cX^{-1}$. These edges form a system of
generators of the groupoid $\cG$, i.e., every morphism $\eta\in\Mor(a,c)$ can be
represented as an admissible composition of generators,
$$
\eta=\xi_{1}*\xi_{2}*\cdots*\xi_{k},\quad \xi_{i}\in\cY\quad 1\leq i\leq k.
$$

Two vertices in the Cayley complex of the groupoid $\cG$ are joined by edges only if
the vertices belong to the same conjugacy class, i.e., when $a,b\in\langle
c\rangle.$Hence, it suffices to consider only the part $\cG_{\langle
c\rangle}$of the groupoid rather than the whole groupoid. Denote the
corresponding Cayley complex by $\cK(\cG_{\langle c\rangle}).$Two edges
$\xi=\left(\frac{a\mapr{}b}{y}\right)$and
$\xi=\left(\frac{b\mapr{}a}{y^{-1}}\right)$are assumed to be equal and have
opposite orientations on the edges.

The one-dimensional cells $\cK_{1}(\cG)$ belong naturally to the set of all morphisms
$\varphi:\cK_{1}(\cG)\hookrightarrow\Mor(\cG)$.

Finally, the two-dimensional cells $\cK_{2}(\cG)$ are the planar orientable polygons
$\sigma(\rho)$given by words $\rho\in\cR\sqcup\cR^{-1}$that define the
boundaries of the polygons $\sigma(\rho)$as closed cycles composed of the
edges of the words $\rho.$The cells $\sigma(\rho)$and
$\sigma(\rho^{-1})$are assumed to be equal and have opposite orientations. The
two-dimensional cells $\sigma(\rho)$are pasted to the 1-skeleton of the
complex $\cK(\cG_{[c]}$by the natural identification of the edges of the
boundary of a cell $\sigma(\rho)$to the corresponding edge of the complex
$\cK(\cG_{[c]}),$preserving the orientation.

\subsection*{Groups of chains of the Cayley complex of the groupoid $\cG$}

The two-dimensional Cayley complex $\K(\cG)$generates the cochain complex $$
C^{0}(\cK(\cG))\mapr{d_{0}}C^{1}(\cK(\cG))\mapr{d_{1}}C^{2}(\cK(\cG)). $$ This
cochain complex has a natural subcomplex of finitely supported cochains, because
every cell of dimension 0 or 1 satisfies the condition that the set of cells that
abut on cells of lesser dimension is finite.

Indeed, if $a\in\Obj(\cG)$is an arbitrary vertex of the Cayley complex
$\cK(\cG)$and $\xi=\left(\frac{a\mapr{}b}{y}\right)$is an arbitrary
edge beginning at $a$, then there are only finitely many edges of this kind,
since $y\in Y,$and the set $Y=X\sqcup X^{-1}$is finite.

Further, if $\xi=\left(\frac{a\mapr{}b}{y}\right)$is an edge, then the cells
of the form of the words $\rho_{i,a_{1}},$
$$ \rho_{i,a_{1}}=
\left(\frac{a_{1}\mapr{}a_{2}}{y_{i1}}\right)
\left(\frac{a_{2}\mapr{}a_{3}}{y_{i2}}\right)
\left(\frac{a_{3}\mapr{}a_{4}}{y_{i3}}\right)\cdots
\left(\frac{a_{l_{i}}\mapr{}a_{1}}{y_{il_{i}}}\right), $$ abut on $\xi$, and,
for some subscripts $i,j$and an element $a_{j},$the condition $y=y_{ij},$$ a_{j}=a$holds. Since the subscripts $i,j$range over a
finite set, it follows that the element $a_{1}$is expressed using $a$in finitely many ways only. Thus, only finitely many words $\rho_{i,a_{1}}$abut on the edge $\xi.$
This, taken together, gives the commutative diagram $$ \xymatrix{
C^{0}(\cK(\cG))\ar[r]^{d_{0}}&C^{1}(\cK(\cG))\ar[r]^{d_{1}}&C^{2}(\cK(\cG))\\
C^{0}_{f}(\cK(\cG))\ar[r]^{d_{0}^{f}}\ar[u]^{\cup}&
C^{1}_{f}(\cK(\cG))\ar[r]^{d_{1}^{f}}\ar[u]^{\cup}& C^{2}_{f}(\cK(\cG)).\ar[u]^{\cup}
} $$ We identify the one-dimensional finitely supported cochains
$C^{1}_{f}(\cK(\cG))$ with the derivations $\Der(C[G])$ by the composition of the
mappings
$$
H:\Der(C[G])\mapr{T}\T_{f}(\cG)\mapr{ \varphi^{*}}C^{1}(\cK(\cG)).
$$
\begin{theorem}\label{t4}
The homomorphism $H$ is a monomorphism onto the kernel of the differential $d_{1}$:
$$
\im(H)=\ker(d_{1}^{f})\subset C^{1}(\cK(\cG)).
$$
The image of the algebra of inner derivations $H(\Int(C[G]))\subset C^{1}(\cK(\cG))$
is equal to the image of the differential $d_{0}^{f}$:
$$
H(\Int(C[G]))=\im(d_{0}^{f})\subset C^{1}(\cK(\cG)).
$$
\end{theorem}

\proof\ 1) To prove that the mapping $H$ is monomorphic, it suffices to show the
monomorphic property of the restriction operator $\varphi^{*}$, since the
correspondence $T$ is an isomorphism (by Theorem~\ref{t3}). Thus, if $T^{X}\in
\T_{f}(\cG)$ and $\varphi^{*}(T^{X})=0$, we are to show that $T^{X}=0$. The character
$T^{X}$ is a function on the groupoid $\cG$, i.e., on $\Mor(\cG)$. The condition
$T^{X}=0$, which is to be proved, means that $T^{X}(\xi)=0$ for every
$\xi\in\Mor(\cG)$. Every morphism $\xi=\left(\frac{a\mapr{}b}{g}\right)$ can be
expanded into an admissible composition $$
\xi=\xi_{1}*\xi_{2}*\xi_{3}*\cdots*\xi_{l}, $$ where the morphisms $\xi_{i}$ are of
the form $$ \xi_{i}=\left(\frac{a_{i}\mapr{}a_{i+1}}{y_{i}}\right),
\quad\xi_{i}\in\cY,\quad 1\leq i\leq l, \quad a=a_{1}, b=a_{l+1}, $$ i.e.,
$\xi_{i}\in \cK_{1}(\cG)$. By assumption,
$T^{X}(\xi_{i})=\varphi^{*}(T^{X})(\xi_{i})=0$. Hence,
$$
T^{X}(\xi)=T^{X}(\xi_{1}*\xi_{2}*\xi_{3}*\cdots*\xi_{l})=\sum\limits_{i=1}^{l}T^{X}(\xi_{i})=0.
$$
2) We claim that $d^{f}_{1}(H(X))=0$ for $X\in\Der(C[G]),$ or
$d^{f}_{1}(\varphi^{*}(T^{X}))=0$ for $T^{X}\in\T_{f}(\cG)$. Strictly speaking, we
are to prove that, on every two-dimensional cell $\sigma(\rho)\in\cK_{2}(\cG)$, the
value of the cochain $d^{f}_{1}(\varphi^{*}(T^{X}))$ vanishes. By the construction of
the two-dimensional Cayley complex, this value is equal to the sum of values of the
cochain $\varphi^{*}(T^{X})$ on the closed cycle formed by the edges of the word
$\rho$. Since the word $\rho$ is one of the words of the set of relations
$\cR\sqcup\cR^{-1}$, $\rho\in\cR\sqcup\cR^{-1}$, we have
$$
\rho=\left(\frac{a_{1}\mapr{}a_{2}}{y_{1}}\right)
\left(\frac{a_{2}\mapr{}a_{3}}{y_{2}}\right)\cdots
\left(\frac{a_{l}\mapr{}a_{1}}{y_{l}}\right),
$$
where the sequence
$$
y_{1}y_{2}\cdots y_{l}\in Y\sqcup Y^{-1}
$$
is one of the defining relations of the group $G$. It follows that
$$
\begin{array}{l}
\varphi^{*}(T^{X})(\rho)=(T^{X})(\rho)=\sum\limits_{i=1}^{l}T^{X}
\left(\frac{a_{i}\mapr{}a_{i+1}}{y_{i}}\right)\\=
T^{X}\left(\left(\frac{a_{1}\mapr{}a_{2}}{y_{1}}\right)*
\left(\frac{a_{2}\mapr{}a_{3}}{y_{2}}\right)*\cdots*
\left(\frac{a_{l}\mapr{}a_{1}}{y_{l}}\right)\right)\\=
T^{X}\left(\frac{a_{1}\mapr{}a_{1}}{y_{1}y_{2}\cdots y_{l}}\right)=
T^{X}\left(\frac{a_{1}\mapr{}a_{1}}{e}\right)=0.
\end{array}
$$
This means that the mapping $\varphi^{*}$ takes the space $\T_{f}(\cG)$ to a subspace
of $\ker(d^{f}_{1})$,
$$
\varphi^{*}: \T_{f}(\cG)\mapr{}\ker(d^{f}_{1}).
$$
3) We claim now that the mapping $\varphi^{*}$ is an epimorphism. The space
$\ker(d^{f}_{1})$ consists of all one-dimensional cocycles of the Cayley complex,
i.e., of the functions on the one-dimensional edges of the Cayley complex that vanish
on every closed one-dimensional chain which is the boundary of a two-dimensional
cell. We are to extend every cochain $\tau$ of this kind to the set $\Mor(\cG)$ to
some character $T$. Let $\xi\in\Mor(\cG)$ be an arbitrary morphism,
$$
\xi=\left(\frac{a\mapr{}b}{g}\right), \quad g\in G.
$$
Every morphism $\xi$ can be expanded in a composition of generating morphisms,
\begin{equation}\label{e4}
\xi=\left(\frac{a\mapr{}a_{2}}{y_{1}}\right)
\left(\frac{a_{2}\mapr{}a_{3}}{y_{2}}\right)\cdots
\left(\frac{a_{l}\mapr{}b}{y_{l}}\right),
\end{equation}
where
$$
y_{1},y_{2},\dots, y_{l}\in X\sqcup X^{-1}.
$$
We write
\begin{equation}\label{e5}
T(\xi)=\sum\limits_{i=1}^{l}\tau\left(\frac{a_{i}\mapr{}a_{i+1}}{y_{i}}\right), \quad
a_{1}=a, \quad a_{l+1}=b.
\end{equation}
The last formula does not depend on the choice of expansion (\ref{e4}), since the
function $\tau$ vanishes on every cycle which is the boundary of a cell. The function
$T$ constructed by formula (\ref{e5}) is obviously a character.

4) Finally, we claim that the image of the algebra of inner derivations
$H(\Int(C[G]))\subset C^{1}(\cK(\cG))$ is equal to the image of the differential
$d_{0}^{f}$. Every inner derivation is defined in the form of a linear combination of
the simplest inner derivations of the form
$$
\ad_{g}:C[G]\mapr{}C[G],\quad \ad_{g}(u)=[g,u], \quad u=\sum\limits_{h\in
G}\lambda^{h}\cdot h\in C[G].
$$
The matrix of the operator $X=\ad_{g}$, $\left\|X_{h}^{h'}\right\|$, is evaluated as
follows:
$$
\sum\limits_{h}\lambda^{h}[g,h]=X(u)=\sum\limits_{h,h'} X_{h}^{h'}\lambda^{h}\cdot
h'.
$$
Since $\lambda^{h}$ are arbitrary, it follows that
$$
[g,h]=gh-hg=\sum\limits_{h'}X^{h'}_{h}\cdot h'.
$$
Then $X^{h'}_{h}=\delta^{h'}_{gh}-\delta^{h'}_{hg}$.

The character $T^{\ad_{g}}$ corresponding to the operator $\ad_{g}$ is a function (on
the groupoid $\cG$) defined by the formula
\begin{equation}\label{e6}
T^{\ad_{g}}\left(\frac{a\mapr{}b}{h}\right)=X^{ha=ah}_{h}=
\delta_{gh}^{ha}-\delta_{hg}^{ha}= \delta_{g}^{b}-\delta_{g}^{a}
\end{equation}

Formula (\ref{e6}) means that
$$
\varphi^{*}\left(T^{\ad_{g}}\right)=d^{f}_{0}(t_{g}), \quad t_{g}\in
C^{0}(\K(\cG)),\quad t_{g}(h)=\delta^{h}_{g}.
$$
This implies that the image of the algebra of inner derivations $H(\Int(C[G]))\subset
C^{1}(\cK(\cG))$ is equal to the image of the differential $d_{0}^{f}$.\qed

\begin{corollary}\label{c1}
The homomorphism $H$ induces an isomorphism of the algebra of outer derivations
$\Out(C[G])$ onto the group of the one-dimensional cohomology with finite supports of
the Cayley complex of the groupoid $\cG$ of the adjoint action of the group~$G$:
$$
H:\Out(C[G])\mapr{\approx} H^{1}_{f}(\cK(\cG); \R).
$$
\end{corollary}

\section{Acknowledgment}
This work was supported by the Russian Foundation for Basic Research (Grant \ No ~ 18-01-00398)

\end{document}